\documentclass[a4,11pt,reqno]{amsart}
\usepackage[english]{babel}    \usepackage[latin1]{inputenc}   \usepackage[T1]{fontenc}     \usepackage[french]{minitoc}
\usepackage[nice]{nicefrac}    \usepackage{latexsym,amsfonts}  
\usepackage{graphics}    \usepackage{ulem}       \usepackage{hhline}    \usepackage{dsfont}    \usepackage{mathrsfs}
\usepackage{fancyhdr}    \usepackage{amsmath}    \usepackage{amssymb}   \usepackage{rotating}  \usepackage{fancybox}
\usepackage{color}       \usepackage{colortbl}   \usepackage{setspace}  \usepackage{enumerate} \usepackage{amsthm}
\usepackage{multicol}    
\usepackage{amsthm}    \usepackage{varioref}  \usepackage{textcomp}
\usepackage{lmodern}     \usepackage{mathpazo}   \usepackage{euscript}  \usepackage[pdftex]{hyperref}
\numberwithin{equation}{section}  \makeatletter\@addtoreset{equation}{section}
   \DeclareMathSymbol{\subsetneqq}{\mathbin}{AMSb}{36}

\newcommand{\C}{\mathbb C}      \newcommand{\R}{\mathbb R}        
           
\begin{document}
\title[]{Coherent states quantization of generalized Bergman spaces on the
unit ball of $\mathbb{C}^{n}$  with a new formula for their associated Berezin transforms}

\author[A. Boussejra, Z. Mouayn]{Abdelhamid Boussejra  \, $\&$ \,  Zouhair Mouayn}
\address{ {(\sf A. B.)}
        {Department of Mathematics,  Faculty of Sciences,
        \newline
         Ibn Tofail University, Kenitra, Morocco
         }}
        \email{boussejra-abdelhamid@univ-ibntofail.ac.ma}

\address{ {(\sf Z. M.)}
        {Department of Mathematics,  Faculty of Sciences and Technics,
        \newline Sultan Moulay Slimane University,BP. 123, Beni Mellal, Morocco
        }}
        \email{mouayn@fstbm.ac.ma}
\date{\today}
\maketitle

\begin{abstract}
While dealing with a class of generalized Bergman spaces on the unit
ball, we construct for each of these spaces a set of coherent states to apply a
coherent states quantization method. This provides us with another way to
recover the Berezin transforms attached to these spaces. Finally, a new
formula representing these transforms a functions of the Laplace-Beltrami
operator is established in terms of Wilson polynomials by using the Fourier-Helgason transform.
\end{abstract}

\section{Introduction}

\bigskip The Berezin transform introduced in \cite{B1} for certain
classical bounded symmetric domains in $\mathbb{C}^{n}$ is a transform linking the
Berezin symbols and symbols for Toeplitz operators. It is present in the
study of the correspondence principle. The formula representing the Berezin
transform as a function of the Laplace operators $\Delta_{1},...,\Delta_{r}$ ( $r$ being the rank of the domain)
plays a key role in the Berezin quantization \cite{B2}.

In this paper, we deal with the rank one symmetric domains. Namely the unit ball $\mathbb{B}^{n}$ in $\left(
\mathbb{C}^{n},\left\langle ,\right\rangle \right) \mathbb{\ }$ endowed with its Bergman metric. We are precisely concerned with
the $L^{2}$-eigenspaces
\begin{equation}
\mathcal{A}_{m}^{2,\nu }\left( \mathbb{B}^{n}\right) =\left\{ \varphi \in
L^{2}(\mathbb{B}^{n},(1-\left| \xi \right| ^{2})^{-n-1}d\mu ),H_{\nu
}\varphi =\epsilon _{m}^{\nu ,n}\varphi \right\}
\end{equation}
associated to the discrete spectrum
\begin{equation}
\epsilon _{m}^{\nu ,n}=4\nu (2m+n)-4m(m+n),m=0,1,2,...,\left[ \nu -n/2\right]
\end{equation}
of the Schr\"{o}dinger operator with uniform magnetic field on $\mathbb{B}%
^{n}$ given by
\begin{equation}
H_{\nu }=-4(1-\left| z\right| ^{2})\left(\sum\limits_{i,j=1}^{n}\left( \delta
_{ij-}z_{i}\overline{z_{j}}\right) \frac{\partial^{2}}{\partial z_{i}\partial\overline{ z}_{j}}+\nu
\sum\limits_{j=1}^{n}(z_{j}\frac{\partial}{\partial z _{j}}-\overline{z}_{j}\frac{\partial}{\partial \overline{z} _{j}})
+\nu ^{2}\right)+4\nu ^{2}
\end{equation}
provided that $\nu >n/2$. Above $\left[ x\right] $ denotes the greatest integer not exceeding
$x.$ For $m\in \mathbb{Z}_{+},$ the Berezin transform associated with the
space in $\left( 1.1\right) $ was obtained in \cite{Gh} via the
well known formalism of Toeplitz operators as
\begin{eqnarray}
\mathfrak{B}_{m}^{\nu ,n}\left[ \varphi \right] \left( z\right)  &=&\frac{%
m!\left( 2\nu -2m-n\right) \Gamma \left( 2\nu -m\right)\Gamma(n) }{n!\Gamma \left( 2\nu -m-n+1\right)\Gamma(n+m) }\int\limits_{\mathbb{B}%
}\left( \frac{(1-\left| z\right| ^{2}(1-\left| \xi \right| ^{2})}{\left|
1-\left\langle z,\xi \right\rangle \right| ^{2}}\right) ^{2\left( \nu
-m\right) }  \notag \\
&&\times \left( P_{m}^{\left( n-1,2(\nu -m)-n\right) }\left( 1-2\left| \xi
\right| ^{2}\right) \right) ^{2}\frac{\varphi \left( \xi \right) }{\left(
1-\left| \xi \right| ^{2}\right) ^{n+1}}d\mu \left( \xi \right)
\end{eqnarray}
where $P_{m}^{\left( \alpha ,\beta \right) }\left( .\right) $ denotes the
Jacobi polynomial \cite{I}. Moreover this transform have been
expressed as a function $f\left( \Delta _{\mathbb{B}^{n}}\right) $ of the
Laplace-Beltrami operator $\Delta _{\mathbb{B}^{n}}$ in terms of an $%
_{3}\digamma _{2}$-sum, see $\left( 5.24\right) $ below.
Our aim here is to construct for each of the eigenspaces in $\left(
1.1\right) $ a set of coherent states by following a generalized formalism \cite{Ga}
in order to apply a coherent states quantization method.
This provides us with another way to recover the Berezin transforms in $%
\left( 1.4\right) $ attached to the $L^{2}$-eigenspace spaces in (1.1).
Finally, we add a new formula expressing the transform $%
\left( 1.4\right) $ as a function of the Laplace-Beltrami operator. The idea
is to make the integral $\left( 1.4\right) $ appear as ''convolution
product'' of the function $\varphi $ with a specific radial function given
in terms of the square of a Jacobi polynomial$.$ Next, a straightforward
computation of the spherical  transform of this radial function with
the use of a Clebsh-Gordon type linearisation \cite{C} for the
square of a Jacobi polynomial\ amounts to a finite sum containing some
integrals whose general form was given by Koornwinder \cite{K} in
terms of Wilson polynomials.

This paper is summarized as follows. In Section 2, we recall briefly the
formalism of coherent states quantization we will be using. Section 3 deals
with some needed facts on the generalized Bergman spaces. In Section 4, we
construct for each of these spaces a set of coherent states and we apply the
corresponding quantization scheme in order to recover their associated
Berezin transforms. In Section 5, we present the formula expressing these
Berezin transforms as functions of the Laplace-Beltrami operator by a different
way and in a new form.

\section{Coherent states quantization}

Coherent states are mathematical tools which provide a close connection
between classical and quantum formalism. In general, they are a specific
overcomplete set of vectors in a Hilbert space satisfying a certain
resolution of the identity condition. Here, we review a coherent states
formalism starting from a measure space ''\textit{as a set of data}'' as
presented in \cite{Ga}. Let $X=\left\{ x\mid x\in X\right\} $ be a
set equipped with a measure $d\mu $ and $L^{2}(X,d\mu )$\ the space of $d\mu
-$square integrable functions on $X.$ Let $\mathcal{A}^{2}$ be a subspace of
$L^{2}(X,d\mu)$  with  an orthonormal basis  $%
\left\{ \Phi _{j}\right\} _{j=0}^{+\infty}$. Let $\mathcal{H}$ be another
(functional) space with  a given orthonormal basis $%
\left\{ \phi _{j}\right\} _{j=0}^{+\infty}$. Then consider the family of states $\left\{ \mid x>\right\}
_{x\in X}$ in $\mathcal{H}$, through the following linear superposition:
\begin{equation}
\mid x>:=\left( \mathcal{N}\left( x\right) \right) ^{-\frac{1}{2}%
}\sum_{j=0}^{+\infty}\Phi _{j}\left( x\right) \mid \phi _{j}>,
\end{equation}
where
\begin{equation}
\mathcal{N}\left( x\right) =\sum_{j=0}^{+\infty}\Phi _{j}\left( x\right) \overline{%
\Phi _{j}\left( x\right) }.
\end{equation}
These coherent states obey the normalization condition
\begin{equation}
\left\langle x\mid x\right\rangle _{\mathcal{H}}=1
\end{equation}
and the following resolution of the identity of $\mathcal{H}$%
\begin{equation}
\mathbf{1}_{\mathcal{H}}=\int\limits_{X}\mid x><x\mid N\left( x\right) d\mu
\left( x\right)
\end{equation}
which is expressed in terms of Dirac's bra-ket notation $\mid x><x\mid $
meaning the rank-one -operator $\varphi \mapsto \left\langle \varphi \mid
x\right\rangle _{\mathcal{H}}.\mid x>.$ The choice of the Hilbert space $%
\mathcal{H}$ define in fact a quantization of the space $X$ by the coherent
states in $\left( 2.1\right) $, via the inclusion map $x\mapsto \mid
x>\in \mathcal{H}$ and the property $\left( 2.4\right) $ is crucial in
setting the bridge between the classical and the quantum mechanics$.$ The
Klauder-Berezin coherent states quantization consists in associating to a
classical observable that is a function $f\left( x\right) $ on $X$ having
specific properties the operator-valued integral
\begin{equation}
A_{f}:=\int\limits_{X}\mid x><x\mid f\left( x\right) \mathcal{N}\left(
x\right) d\mu \left( x\right)
\end{equation}
The function $f\left( x\right) \equiv \widehat{A}_{f}\left( x\right) $ is
called upper (or contravariant) symbol of the operator $A_{f}$ and is
nonunique in general. On the other hand, the expectation value $\left\langle
x\mid A_{f}\mid x\right\rangle $ of $A_{f}$ with respect to the set of
coherent states $\left\{ \mid x>\right\} _{x\in X}$ is called lower ( or
covariant) symbol of $A_{f}.$ Finally, associating to the classical
observable $f\left( x\right) $ the obtained mean value $\left\langle x\mid
A_{f}\mid x\right\rangle ,$ we get the Berezin transform of this observable.
That is,
\begin{equation}
B\left[ f\right] \left( x\right) :=\left\langle x\mid A_{f}\mid
x\right\rangle ,\text{ }x\in X.
\end{equation}
For all aspect of the theory of coherent states and their genesis, we refer
to the survey \cite{D} by Dodonov or to the book by Gazeau \cite{Ga}.

\section{The spaces $\mathcal{A}_{m}^{2,\protect\nu }\left( \mathbb{B}%
^{n}\right) $}

In this section, we review some results on the $L^{2}$-concrete spectral
analysis of the Schr\"{o}dinger operator $H_{\nu }$ in $\left( 1.3\right) $
and acting in  the Hilbert space $L^{2}(\mathbb{B}^{n},d\mu_{n})$,see \cite
{B}, for more details.
Let $\mathbb{B}^{n}=\{z\in \mathbb{C}^{n};\mid z\mid <1\}$ be the unit ball in ${\C}^{n}$ with  the Lebesgue measure $d\mu$ normalized so that $\mu(\mathbb{B}^{n})$  and let $\partial{\mathbb{B}^{n}}=\{\omega\in \mathbb{C}^{n},\mid \omega\mid=1\}$ be the unit sphere with $d\sigma$ the normalized measure on it.
Let $G=SU(n,1)$ be the group of all $\mathbb{C}$-linear
transforms $g$ on $\mathbb{C}^{n+1}$ that preserve the indefinite hermitian
form $\sum_{j=1}^{n}\mid z_{j}\mid ^{2}-\mid z_{n+1}\mid ^{2},$ with $\det
g=1$. Then $G$  acts transitively on the unit ball by
\begin{equation}
G\ni g=\left(
\begin{matrix}
a & b \\
c & d
\end{matrix}
\right):z\rightarrow g.z=(az+b)(cz+d)^{-1}.
\end{equation}
As a homogeneous space we have the identification $\mathbb{B}^{n}=G/K$ where
$K=S(U(n)\times U(1))$ is the stabilizer of $0$.
It is endowed with its
usual Khaler-Bergman metric $ds^{2}=-\sum_{i,j}^{n}\partial _{j}\overline{%
\partial }_{j}(Log(1-\left| z\right| ^{2}))dz_{i}\otimes \overline{dz_{j}}$.
The Bergman distance and the volume element on $\mathbb{B}^{n}$ are given
respectively by
\begin{equation}
\cosh ^{2}d\left( z,w\right) =\frac{\left| 1-\left\langle z,w\right\rangle
\right| ^{2}}{(1-\mid z\mid^{2} )(1-\mid w\mid^{2} )}
\end{equation}
and $d\mu _{n}\left( z\right) =(1-\mid z\mid^{2} )^{-\left(
n+1\right) }d\mu \left( z\right)$.\\
The group $G$ acts  unitarily  on  the space $L^{2}(\mathbb{B}^{n},d\mu_{n})$, via
$
U(g)F(z)=F(g^{-1}.z).
$
Let consider the magnetic gauge vector potential given through the canonical
1-form on $\mathbb{B}^{n}$: $\theta =-i(\partial -\overline{\partial }%
)Log(1-\mid z\mid^{2} )$, to which the Schrodinger operator
\begin{equation}
H_{\nu }=-\left( d+i\nu \text{ext}\left( \theta \right) \right) ^{\ast
}\left( d+i\text{ext}\left( \theta \right) \right) +4\nu ^{2}
\end{equation}
can be associated. Here $\nu \geq 0$ is a fixed number $d$ \ denotes the
usual exterior derivative on differential forms on $\mathbb{B}^{n}$ and ext$%
\left( \theta \right) $ is the exterior multiplication by $\theta $ while
the symbol $\ast $ stands for the adjoint operation with respect to the
Hermitian scalar product induced by the Bergman metric $ds^{2}$ on
differential forms. Note that when $\nu =0,$ the operator in $\left(
3.3\right) $ reduces to
\begin{equation}
H_{0}\equiv \Delta _{\mathbb{B}^{n}}=4(1-\mid z\mid^{2})
 \sum_{i,j=1}^{n}\left( \delta _{ij}-z_{i}\overline{z_{j}}\right) \frac{\partial^{2}}{\partial z_{i}\partial\overline{ z}_{j}}
\end{equation}
which is the Laplace-Beltrami operator of the Bergman ball $\mathbb{B}^{n}$. For general $\nu
\geq 0,$ the Schrodinger operator $H_{\nu }$ in $\left( 3.3\right) $ can be
expressed in the complex coordinates $\left( z_{1},...,z_{n}\right) $ by the
formula $\left( 1.3\right) $ see \cite{Ay},\cite{Bo} and \cite{Ge}.

Now, for an arbitrary complex number $\lambda $, a fundamental family of
eigenfunctions of $H_{\nu }$ with eigenvalue $\lambda ^{2}+4\nu ^{2}+n^{2}$
is given by the Poisson kernels :
\begin{equation}
z\mapsto P_{\lambda }^{\nu }(z,\theta )=\left( \frac{1-\mid z\mid ^{2}}{\mid
1-<z,\theta >\mid ^{2}}\right) ^{\frac{1}{2}\left( i\lambda +1\right)
}\left( \frac{1-\overline{<z,\theta >}}{1-<z,\theta >}\right) ^{\nu },z\in
\mathbb{B}^{n}.
\end{equation}
Moreover, a complete description of the
expansion of an eigenfunction $f$ of $H_{\nu }$ with eigenvalue $\lambda
^{2}+4\nu ^{2}+n^{2}$, in terms of the appropriate Fourier series in $%
\mathbb{B}^{n}$ have been given in \cite[Proposition 2.2]{B}. Precisely,
\begin{equation}\label{e:barwq}
f(z)=(1-\rho ^{2})^{\frac{i\lambda +n}{2}}
\sum\limits_{p,q=0}^{+\infty }\rho ^{p+q} \cdot_{2}\digamma _{1}\left( \frac{i\lambda +n}{2}+\nu +p,\frac{%
i\lambda +n}{2}-\nu +q,p+q+n;\rho ^{2}\right) a_{p,q}^{\lambda ,\nu
}.h_{p,q}(\theta ),
\end{equation}
in $C^{\infty }([0,1[\times\partial{\mathbb{B}^{n}})$, $z=\rho \theta $, $\rho \in \lbrack 0,1[$
and $\mid \theta \mid =1$. Above $_{2}F_{1}$ denotes the Gauss
hypergeometric function \cite{Gr} and $a_{p,q}^{\lambda
,\nu }=(a_{p,q,j}^{\lambda ,\nu })\in {\C}^{d(n,p,q)}$ are complex numbers,
where
\begin{equation}
d(n,p,q):=\frac{(p+q+n-1)(p+n-2)!(q+n-2)!}{p!q!(n-1)!(n-2)!}
\end{equation}
is the dimension of the space $H(p,q)$ of restrictions to the unit sphere $%
\partial \mathbb{B}^{n}$ of harmonic polynomials $h(z)$ on ${\C}^{n}$, which
are homogeneous of degree $p$ in $z$ and degree $q$ in $\overline{z}$, see \cite{F} or \cite{R} for more details.
The notation ''.'' in (3.6) means the following finite sum
\begin{equation}
a_{p,q}^{\lambda ,\nu }.h_{p,q}(\theta
)=\sum\limits_{j=1}^{d(n,p,q)}a_{p,q,j}^{\lambda ,\nu }h_{p,q}^{j}(\theta ),
\end{equation}
where $\{h_{p,q}^{j}\}_{1\leq j\leq d(n,p,q)}$ is an orthonormal basis of
$H(p,q)$. The spectral analysis of $H_{\nu }$ have been studied by many
authors, see \cite{B} and references therein.
Actually, $H_{\nu }$ is an elliptic densely defined operator on the Hilbert
space $L^{2}(\mathbb{B}^{n},(1-\left\langle z,z\right\rangle )^{-\left(
n+1\right) }d\mu )$ admitting a unique self-adjoint realization also denoted
by $H_{\nu }$. Its spectrum consists of a continuous part given by $\left[
n^{2},+\infty \right[ $ (corresponding to scattering states) and a finite
number of infinitely degenerate eigenvalues $\epsilon _{m}^{\nu ,n}$ given
by $\left( 1.2\right) $ (characterizing bound states) provided that $2\nu >n.
$ More precisely, $\epsilon _{m}^{\nu ,n}=\lambda _{m}^{2}+4\nu ^{2}+n^{2}$,
with $\lambda _{m}=i(2m+n-2\nu )$, $m=0,1,...,\left[ \nu -n/2\right] $.
Here, we focus on the discrete part of the spectrum, which is labeled by
the integer $m$ and the corresponding eigenspace $\mathcal{A}_{m}^{2,\nu
}\left( \mathbb{B}^{n}\right) $ defined in $\left( 1.1\right) $. Taking into
account (3.6) and expressing the involved hypergeometric in terms of Jacobi
polynomial, an orthonormal basis of $\mathcal{A}_{m}^{2,\nu }\left( \mathbb{B%
}^{n}\right) $ can be given explicitly by
\begin{equation}
\Phi _{p,q}^{\nu ,m,j}\left( z\right) =\kappa _{p,q}^{\nu ,m,n}\left(
1-\left| z\right| ^{2}\right) ^{\nu -m}P_{m-q}^{( n+p+q-1,2(\nu -m)
-n) }\left( 1-2\left| z\right| ^{2}\right) h_{p,q}^{j}\left( z,%
\overline{z}\right)
\end{equation}
with
\begin{equation}
\kappa _{p,q}^{\nu ,m,n}=\left(\frac{n\Gamma(2\nu-m-n-q+1)\Gamma(p+n+m)}
{(m-q)!(2(\nu-m)-n)\Gamma(2\nu-m+p)}\right)^{-\frac{1}{2}}.
  \end{equation}
for varying $p=0,1,2,...$, $q=0,1,....,m$ and $j=1,....,d(n;p,q)$.
Furthermore, the space $\mathcal{A}_{m}^{2,\nu }(\mathbb{B}^{n}) $ is a reproducing kernel Hilbert space. That is, there exists a unique complex valued function $K^{\nu,m}$ on $\mathbb{B}^{n}\times\mathbb{B}^{n}$ such that, denoting $K^{\nu,m}_{z}(w)=K^{\nu,m}(w,z)$, $K^{\nu,m}_{z}$ belongs to $\mathcal{A}_{m}^{2,\nu }(\mathbb{B}^{n})$ for any $z\in\mathbb{B}^{n}$ and
$$
f(z)=<f,K^{\nu,m}_{z}>,
$$
for all functions $f$ in $\mathcal{A}_{m}^{2,\nu }(\mathbb{B}^{n})$ and all $z\in \mathbb{B}^{n}$.Its expression can be given explicitly as function of the Bergman geodesic distance as
\begin{eqnarray}
K^{\nu,m}(z,w)  &=&%
\frac{\left( 2\left( \nu -m\right) -n\right) \Gamma
\left( 2\nu -m\right) }{n! \Gamma \left( 2\nu
-m-n+1\right) }\left( \frac{(1-\overline{< z,w>})}{1-<z,w>}\right) ^{\nu }  \\
&&\times (\cosh d(z,w))^{-2(\nu-m)})  P_{m}^{\left( n-1,2(\nu -m)-n\right) }( 1-2\tanh^{2}d(z,w)) \notag
\end{eqnarray}

\noindent \textbf{Remark 3.1}. For $m=0,$ the space $\mathcal{A}_{0}^{2,\nu }\left(
\mathbb{B}^{n}\right) $ reduces further to be isomorphic to the weighted
Bergman space of holomorphic function $\psi $ on $\mathbb{B}^{n}$ satisfying
the growth condition
\begin{equation*}
\int_{\mathbb{B}^{n}}\left| \psi \left( z\right) \right|
^{2}((1-\left\langle z,z\right\rangle )^{2\nu -n-1}d\mu \left( z\right)
<+\infty .
\end{equation*}
This fact justify why the eigenspace $\mathcal{A}_{m}^{2,\nu }\left( \mathbb{%
B}^{n}\right) $ have been also called a\textit{\ generalized Bergman spaces
of index }$m$.

\section{Coherent states quantization}

Now, to adapt the defintion (2.1) of coherent states for the context of the
generalized Bergman spaces in $\left( 1.1\right) $ we first list the
following notations.
\begin{itemize}
\item  $\left( X,d\eta \right) :=\left( \mathbb{B}^{n},\left( 1-\left|
z\right| ^{2}\right) ^{-\left( n+1\right) }d\mu \right) ,d\eta \equiv d\mu
_{n}$ is the volume element on $\mathbb{B}^{n}.$

\item  $x\equiv z\in \mathbb{B}^{n}.$

\item  $\mathcal{A}^{2}:=\mathcal{A}_{m}^{2,\nu }\left( \mathbb{B}%
^{n}\right) \subset L^{2}(\mathbb{B}^{n},\left( 1-\left| z\right|
^{2}\right) ^{-n-1}d\mu ).$

\item  $\left\{ \Phi _{k}\left( x\right) \right\} \equiv \left\{ \Phi
_{p,q,j}^{\nu ,m}\left( z\right) \right\} $ is the orthonormal basis of $%
\mathcal{A}_{m}^{2,\nu }\left( \mathbb{B}^{n}\right) $ in $\left( 3.8\right)
$

\item  $\mathcal{N}\left( x\right) \equiv \mathcal{N}\left( z\right) $ is a
normalization factor.

\item  $\left\{ \varphi _{k}\right\} \equiv \left\{ \varphi _{p,q,j}\right\} $
is an orthonormal basis of another (functional) Hilbert space $\mathcal{H}$.
\end{itemize}

\quad

\noindent \textbf{Definition 4.1.} \textit{For each fixed integer }$m=0,1,...,\left[
n-\nu /2\right] ,$\textit{\ a class of generalized coherent states associated
with the space }$A_{m}^{2,\nu }\left( \mathbb{B}^{n}\right) $\textit{\ is
defined according to }$\left( 2.1\right) $\textit{\ by the form}
\begin{equation}
\phi _{z}^{\nu ,m}\equiv \mid z,\nu ,m>:=\left( \mathcal{N}\left( z\right)
\right) ^{-\frac{1}{2}}\sum\limits_{\substack{ 0\leq q\leq m,0\leq p<+\infty
\\ 1\leq j\leq d\left( n,p,q\right) }}\Phi _{p,q,j}^{\nu ,m}\left( z\right)
\varphi _{p,q,j}
\end{equation}
\textit{where }$\mathcal{N}\left( z\right) $\textit{\ is a normalization
factor.}
\\

\noindent \textbf{Proposition 4.1. }\textit{The factor in }$\left( 4.1\right) $\textit{%
\ is given by }
\begin{equation}
\mathcal{N}\left( z\right) =\frac{\left( 2(\nu -m)-n\right) \Gamma \left(
2\nu -m\right) }{n!\Gamma \left( 2\nu -m-n+1\right) }\frac{\Gamma
\left( m+n\right) }{m!\Gamma \left( n\right) }
\end{equation}
\textit{for every} $z\in \mathbb{B}^{n}.$
\\

 \noindent \textbf{Proof.} To calculate this factor, we start by writing the
condition
\begin{equation}
\left\langle \phi _{z}^{\nu ,m},\phi _{z}^{\nu ,m}\right\rangle _{\mathcal{H}%
}=1.
\end{equation}
Equation $\left( 4.3\right) $ is equivalent to
\begin{equation}
\left( \mathcal{N}\left( z\right) \right) ^{-1}\sum\limits_{p=0}^{+\infty
}\sum\limits_{q=0}^{m}\sum\limits_{j=1}^{d(n,p,q)}\Phi _{p,q,j}^{\nu ,m}\left( z\right) \overline{\Phi
_{p,q,j}^{\nu ,m}\left( z\right) }=1
\end{equation}
Making use of (3.9) and (3.11) for the particular case $z=w,$ we get that
\begin{equation}
\mathcal{N}\left( z\right) =\frac{\left( 2(\nu -m)-n\right) \Gamma \left(
2\nu -m\right) }{n!\Gamma \left( 2\nu -m-n+1\right) }P_{m}^{\left(
n-1,2(\nu -m)-n\right) }\left( 1\right)
\end{equation}
Next, by the following fact on Jacobi polynomial \cite{Gr}:
\begin{equation}
P_{m}^{\left( \alpha ,\beta \right) }\left( 1\right) =\frac{\Gamma \left(
m+\alpha +1\right) }{m!\Gamma \left( \alpha +1\right) }
\end{equation}
for $\alpha =n-1$ to arrive at the announced result.The states $\phi _{z}^{\nu ,m}\equiv \mid z,\nu ,m>$ satisfy the
resolution of the identity
\begin{equation}
1_{\mathcal{H}}=\int\limits_{\mathbb{B}^{n}}\mid z,\nu ,m><z,\nu ,m\mid
\mathcal{N}\left( z\right) d\nu
\end{equation}
and with the help of them we can achieve the coherent states quantization
scheme described in Sec.2 to rederive the Berezin transform $\frak{B}%
_{m}^{\nu ,n}$ in $\left( 1.4\right) $ which was defined by Toeplitz
operators formalism in \cite{Gh}. For this let us associate to any
arbitrary function $\varphi \in L^{2}(\mathbb{B}^{n},(1-\left| \xi \right|
^{2})^{-n-1}d\mu )$ the operator-valued integral
\begin{equation}
A_{\varphi }:=\int\limits_{\mathbb{B}^{n}}\mid z,\nu ,m><z,\nu ,m\mid
\varphi \left( z\right) \mathcal{N}\left( z\right) (1-\left| z\right|
^{2})^{-n-1}d\mu
\end{equation}
The function $\varphi \left( z\right) $ is a upper symbol of the operator $%
A_{\varphi }.$ On the other hand, we need to calculate the expectation value
\begin{equation}
\mathbb{E}_{\left\{ \mid z,\nu ,m>\right\} }\left( A_{\varphi }\right)
:=<z,\nu ,m\mid A_{\varphi }\mid z,\nu ,m>
\end{equation}
of \ $A_{\varphi }$ with respect to the set of coherent states $\left\{ \mid
z,\nu ,m>\right\} _{z\in \mathbb{B}^{n}}$ defined in $\left( 4.1\right)$.
This will constitute a lower symbol of the operator \ $A_{\varphi }.$\\

\noindent \textbf{Proposition 4.2.} \textit{Let }$\varphi \in L^{2}(\mathbb{B}%
^{n},(1-\left| \xi \right| ^{2})^{-n-1}d\mu ).$\textit{\ Then, the
expectation value in }$\left( 4.9\right) $\textit{\ has the following
expression }
\begin{eqnarray}
\mathbb{E}_{\left\{ \mid z,\nu ,m>\right\} }\left( A_{\varphi }\right)  &=&%
\frac{\Gamma \left( n\right) m!\left( 2\left( \nu -m\right) -n\right) \Gamma
\left( 2\nu -m\right) }{n!\Gamma \left( n+m\right) \Gamma \left( 2\nu
-m-n+1\right) }\int\limits_{\mathbb{B}}\left( \frac{(1-\left| z\right|
^{2}(1-\left| \xi \right| ^{2})}{\left| 1-\left\langle z,\xi \right\rangle
\right| ^{2}}\right) ^{2\left( \nu -m\right) }  \\
&&\times \left( P_{m}^{\left( n-1,2(\nu -m)-n\right) }\left( 1-2\left| \xi
\right| ^{2}\right) \right) ^{2}\frac{\varphi \left( \xi \right) }{\left(
1-\left| \xi \right| ^{2}\right) ^{n+1}}d\mu \left( \xi \right)   \notag
\end{eqnarray}
\textit{for every }$z\in \mathbb{B}^{n}.$\\

\noindent \textbf{Proof. }We first write the action of the operator $A_{\varphi }$ in $%
\left( 4.8\right) $ on an arbitrary coherent state $\mid z,\nu ,m>$ in terms
of Dirac's bra-ket notation as
\begin{equation}
A_{\varphi }\mid z,\nu ,m>=\int\limits_{\mathbb{B}^{n}}\mid w,\nu ,m><w,\nu
,m\mid z,\nu ,m>\frac{\mathcal{N}\left( w\right) }{(1-\left| w\right|
^{2})^{n+1}}d\mu \left( w\right)
\end{equation}
Therefore, the expectation value reads
\begin{align}
 <z,\nu ,m\mid A_{\varphi }\mid z,\nu ,m>
& \qquad =\int\limits_{\mathbb{B}^{n}}<z,\nu ,m\mid w,\nu ,m>\overline{<z,\nu ,m\mid
w,\nu ,m>}\frac{\mathcal{N}\left( w\right) }{(1-\left| w\right| ^{2})^{n+1}}%
d\mu \left( w\right)
\\
& \qquad =\int\limits_{\mathbb{B}^{n}}\left| <z,\nu ,m\mid w,\nu ,m>\right|
^{2}\varphi \left( w\right) \frac{\mathcal{N}\left( w\right) }{(1-\left|
w\right| ^{2})^{n+1}}d\mu \left( w\right).
\end{align}
Now, we need to evaluate the quantity $\left| <z,\nu ,m\mid w,\nu ,m>\right|
^{2}$ in $\left( 4.13\right)$. For this, we write the scalar product as
\begin{equation}
<z,\nu ,m\mid w,\nu ,m>=\sum_{p=0}^{+\infty
}\sum_{q=0}^{m}\sum\limits_{j=1}^{d\left( n;p,q\right) }\sum_{r=0}^{+\infty
}\sum_{s=0}^{m}\sum\limits_{l=1}^{d\left( n;p,q\right) }\frac{\Phi
_{p,q,j}^{\nu ,m}\left( z\right) \overline{\Phi _{p,q,j}^{\nu ,m}\left(
w\right) }}{\sqrt{\mathcal{N}\left( z\right) \mathcal{N}\left( w\right) }}%
\left\langle \varphi _{p,q,j},\varphi _{p,q,l}\right\rangle _{\mathcal{H}}
\end{equation}
Recalling that
\begin{equation}
\left\langle \varphi _{p,q,j},\varphi _{p,q,l}\right\rangle _{\mathcal{H}%
}=\delta _{j,l}\delta _{p,r}\delta _{q,s}
\end{equation}
since $\left\{ \varphi _{p,q,j}\right\} $ is an orthonormal basis of $%
\mathcal{H}$, the above sum in $\left( 4.14\right) $ reduces to
\begin{equation}
<z,\nu ,m\mid w,\nu ,m>=\left( \mathcal{N}\left( z\right) \mathcal{N}\left(
w\right) \right) ^{-\frac{1}{2}}\sum\limits_{\substack{ 0\leq q\leq m,0\leq
p<+\infty  \\ 1\leq j\leq d\left( n,p,q\right) }}\Phi _{p,q,j}^{\nu
,m}\left( z\right) \overline{\Phi _{p,q,j}^{\nu ,m}\left( w\right) }.
\end{equation}
Now, taking account (3.9) and (3.11), Equation ( 4.16) takes the form
\begin{align}
<z,\nu ,m\mid w,\nu ,m> &=\frac{\left( 2\left( \nu -m\right) -n\right) \Gamma
\left( 2\nu -m\right) }{n!\Gamma \left( 2\nu -m-n+1\right) }\left(
\mathcal{N}\left( z\right) \mathcal{N}\left( w\right) \right) ^{-\frac{1}{2}%
}\left( \frac{1-\overline{\left\langle z,w\right\rangle }}{1-\left\langle
z,w\right\rangle }\right) ^{\nu }
\\
&
\times \left( \cosh \left( d\left( z,w\right) \right) \right) ^{-2\left( \nu
-m\right) }P_{m}^{\left( n-1,2\left( \nu -m\right) -n\right) }\left(
1-2\tanh ^{2}\left( d\left( z,w\right) \right) \right). \nonumber
\end{align}
So that the square modulus of the scalar product in $\left( 4.17\right) $
reads
\begin{align}
\left| <z,\nu ,m\mid w,\nu ,m>\right| ^{2}
&=\left( \frac{\left( 2\left( \nu -m%
\right) -n\right) \Gamma \left( 2\nu -m\right) }{n!\Gamma \left( 2\nu
-m-n+1\right) }\right) ^{2}\left( \mathcal{N}\left( z\right) \mathcal{N}%
\left( w\right) \right) ^{-1}
\\
&
\times \left( \cosh \left( d\left( z,w\right) \right) \right) ^{-4\left( \nu
-m\right) }\left( P_{m}^{\left( n-1,2\left(\nu -m\right) -n\right) }\left(
1-2\tanh ^{2}\left( d\left( z,w\right) \right) \right) \right) ^{2}. \nonumber
\end{align}
Returning back to $\left( 4.12\right),$ we get
\begin{align}
\mathbb{E}_{\left\{ \mid z,\nu ,m>\right\} }\left( A_{\varphi }\right)
&=\int\limits_{\mathbb{B}^{n}}\varphi \left( w\right) \left( \frac{\left( 2%
\left[ \nu -m\right] -n\right) \Gamma \left( 2\nu -m\right) }{n!\Gamma
\left( 2\nu -m-n+1\right) }\right) ^{2}\left( \mathcal{N}\left( z\right)
\mathcal{N}\left( w\right) \right) ^{-1}\frac{\mathcal{N}\left( w\right) }{%
(1-\left| w\right| ^{2})^{n+1}}
\\&
\times \left( \cosh \left( d\left( z,w\right) \right) \right) ^{-4\left( \nu
-m\right) }\left( P_{m}^{\left( n-1,2\left( \nu -m\right) -n\right) }\left(
1-2\tanh ^{2}\left( d\left( z,w\right) \right) \right) \right) ^{2}d\mu
\left( w\right), \nonumber
\end{align}
which can be also written as
\begin{align}
\mathbb{E}_{\left\{ \mid z,\nu ,m>\right\} }\left( A_{\varphi }\right)
&=\int\limits_{\mathbb{B}^{n}}\varphi \left( w\right) \left( \frac{\left( 2%
\left( \nu -m\right) -n\right) \Gamma \left( 2\nu -m\right) }{n!\Gamma
\left( 2\nu -m-n+1\right) }\right) ^{2}\frac{\left( \mathcal{N}\left(
z\right) \right) ^{-1}}{(1-\left| w\right| ^{2})^{n+1}}
\\
&\qquad \times \left( \cosh \left( d\left( z,w\right) \right) \right) ^{-4\left( \nu
-m\right) }\left( P_{m}^{\left( n-1,2\left( \nu -m\right) -n\right) }\left(
1-2\tanh ^{2}\left( d\left( z,w\right) \right) \right) \right) ^{2}d\mu
\left( w\right)
\nonumber \\
& =\frac{\left( 2(\nu -m)-n\right) \Gamma \left( 2\nu -m\right) m!\Gamma \left(
n\right) }{n!\Gamma \left( 2\nu -m-n+1\right) \Gamma \left( m+n\right)
}\int\limits_{\mathbb{B}^{n}}\frac{%
\varphi \left( w\right) }{(1-\left| w\right| ^{2})^{n+1}}
\\
& \qquad \times \left( \cosh \left( d\left( z,w\right) \right) \right) ^{-4\left( \nu
-m\right) }\left( P_{m}^{\left( n-1,2\left( \nu -m\right) -n\right) }\left(
1-2\tanh ^{2}\left( d\left( z,w\right) \right) \right) \right) ^{2}d\mu
\left( w\right). \nonumber
\end{align}
Finally, we summarize the above discussion by considering the following
definition.\\

\noindent \textbf{Definition 4.3.} \textit{The Berezin transform of the classical
observable }$\varphi \in L^{2}(\mathbb{B}^{n},(1-\left| \xi \right|
^{2})^{-n-1}d\mu )$\textit{\ constructed according to the quantization by
the coherent states} $\left\{ \mid z,\nu ,m>\right\} $ \textit{in }$\left(
4.1\right) $\textit{\ is obtained \ by} \textit{associating to }$\varphi $%
\textit{\ the obtained mean value in }$\left( 4.10\right) .$\textit{\ That
is, }
\begin{equation}
\frak{B}_{m}^{\nu ,n}\left[ \varphi \right] \left( z\right) =\mathbb{E}%
_{\left\{ \mid z,\nu ,m>\right\} }\left( A_{\varphi }\right)
\end{equation}
\textit{for every }$z\in \mathbb{B}^{n}.$\\

\noindent \textbf{Remark 4.4}. For $m=0,$ the transform $\left( 4.10\right) $ reduces
to the well known Berezin transform attached to the weighted Bergman space $%
\mathcal{A}_{0}^{2,\nu }\left( \mathbb{B}^{n}\right) $ of holomorphic
function $\psi $ on $\mathbb{B}^{n}$ satisfying the growth condition $\left(
3.12\right) $ and given by
\begin{equation}
\frak{B}_{0}^{\nu ,n}\left[ \varphi \right] \left( z\right) =\frac{\left(
2\nu -n\right) \Gamma \left( 2\nu \right) }{n!\Gamma \left( 2\nu
-n+1\right) }\int\limits_{\mathbb{B}^{n}}\left( \cosh d\left( z,\xi \right)
\right) ^{-4\nu }\frac{\varphi \left( \xi \right) }{\left( 1-\left| \xi
\right| ^{2}\right) ^{n+1}}d\mu \left( \xi \right)
\end{equation}
The latter one have also been written as a function of the Bergman Laplacian
$\Delta _{\mathbb{B}^{n}}$ as
\begin{equation}
\frak{B}_{0}^{\nu ,n}=\frac{1}{\Gamma \left( \alpha +1\right) \Gamma \left(
\alpha +n+1\right) }\left| \Gamma \left( \alpha +1+\frac{n}{2}+\frac{i}{2}%
\sqrt{-\Delta _{\mathbb{B}^{n}}-n^{2}}\right) \right| ^{2}
\end{equation}
firstly by Berezin. The above form,
involving gamma factors, was derived by Peetre in \cite[p. 182]{P},  so that $\alpha $ there occurring in the weight of the Bergman space, corresponds to $2\nu -n-1.$

\section{An expression of $\frak{B}_{m}^{\protect\nu ,n}$ as function of $%
\Delta _{\mathbb{B}^{n}}$}
Then Berezin transform $\frak{B}_{m}^{\nu ,n}$ associated the generalized Bergman space  $\mathcal{A}_{m}^{2,\nu
}\left( \mathbb{B}^{n}\right)$ is given  by
\begin{equation}
\frak{B}_{m}^{\nu ,n}\left[ \varphi \right] \left( z\right) =c_{m}^{\nu
,n}\int\limits_{\mathbb{B}^{n}}\frac{\left( P_{m}^{\left( n-1,2(\nu
-m)-n\right) }\left( 1-2\tanh ^{2}d\left( z,\xi \right) \right) \right)^{2}%
}{\left( \cosh d\left( z,\xi \right) \right) ^{4\left( \nu -m\right) }}%
\varphi \left( \xi \right) \frac{d\mu \left( \xi \right) }{\left( 1-\left|
\xi \right| ^{2}\right) ^{n+1}},
\end{equation}
with
\begin{equation}
c_{m}^{\nu ,n}=\frac{\Gamma \left( n\right) m!\left( 2\left( \nu -m\right)
-n\right) \Gamma \left( 2\nu -m\right) }{n!\Gamma \left( n+m\right)
\Gamma \left( 2\nu -m-n+1\right) }
\end{equation}
Let $B^{\nu,n}_{m}(z,w)$ be the kernel function of the above integral operator and set  $h^{\nu,n}_{m}(g)=B^{\nu,n}_{m}(z,0)$, $z=g.0$.
Then the integral operator (5.1) can be written as a convolution product over $G$:
$$
\frak{B}_{m}^{\nu ,n}\left[ \varphi \right] \left( z\right)= c_{m}^{\nu,n}(\varphi\ast h^{\nu,n}_{m})(g),\quad z=g.0,
$$
from which it follows easily that the Berezin operator is an  $L^{2}$-bounded operator.\\
Since  $B^{\nu,n}_{m}(z,w)$ is a $G$ bi-invariant function it follows that  $\frak{B}_{m}^{\nu ,n}$ is a $G$-invariant operator. That is $U(g)\circ \frak{B}_{m}^{\nu ,n}=\frak{B}_{m}^{\nu ,n}\circ U(g)$, for every $g\in G$. Therefore  $\frak{B}_{m}^{\nu ,n}$  is, in the spectral theoretic sense,  a function of the $G$-invariant Laplacian  $\Delta _{\mathbb{B}^{n}}$ of the unit ball. Below we give it explicitly.\\

\noindent \textbf{Theorem 5.1.}{\it The Berezin transform $\frak{B}_{m}^{\nu ,n}$
can be expressed as a function of the Laplace-Beltrami operator $%
\Delta _{\mathbb{B}^{n}\text{ }}$ as
\begin{equation}\label{e:barwq}\begin{split}
\frak{B}_{m}^{\nu ,n}&=\left| \Gamma \left( 2\left( \nu -m\right) -\frac{n%
}{2}+\frac{i}{2}\sqrt{-\Delta _{\mathbb{B}^{n}}-n^{2}}\right) \right| ^{2}\\
&\quad \sum\limits_{k=0}^{2m}\gamma_{k}^{\nu,m,n}W_{k}(-\frac{1}{4}\Delta _{\mathbb{B}^{n}}-\frac{n^2}{4};2(\nu-m)-\frac{n}{2},1+\frac{n}{2},\frac{n}{2},\frac{n}{2})
\end{split}\end{equation}
where $W_{k}(.)$ are Wilson polynomials,
$$
 \gamma _{k}^{\nu ,n,m}=\frac{2m!\Gamma \left( n\right) \left( 2\left( \nu
-m\right) -n\right) \Gamma \left( 2\nu -m\right) \left( -1\right) ^{k}}{
\Gamma \left( n+m\right) \Gamma \left( 2\nu -m-n+1\right) k!\Gamma
^{2}\left( 2\left( \nu -m\right) +k\right) }\times A_{k}^{\nu ,n,m},
$$
and the coefficients $A_{k}^{\nu ,n,m}$ are given by (5.10) below.} \\

\noindent \textbf{Proof. } Recall that
$$
\frak{B}_{m}^{\nu ,n}\left[ \varphi \right] \left( z\right)= c_{m}^{\nu,n}(\varphi\ast h^{\nu,n}_{m})(g),\quad z=g.0,
$$
where
\begin{equation}
h_{m}^{\nu ,n}\left( \xi \right) :=\left( 1-\left| \xi \right| ^{2}\right)
^{2\left( \nu -m\right) }\left( P_{m}^{\left( n-1,2(\nu -m)-n\right) }\left(
1-2\left| \xi \right| ^{2}\right) \right) ^{2},\xi \in \mathbb{B}^{n},
\end{equation}
By this way, we have to compute  the spherical transform $\mathcal{F}[h_{m}^{\nu ,n}]$ of $h_{m}^{\nu ,n}$. Namely
\begin{equation}
\mathcal{F}[h_{m}^{\nu ,n}](\lambda):=\int\limits_{\mathbb{B}^{n}}h_{m}^{\nu ,n}(z)
\phi_{-\lambda}(z)d\mu_{n}(z),\lambda\in {\R}
\end{equation}
where $\phi_{\lambda}$ is the spherical function associated to $\Delta _{\mathbb{B}^{n}}$, given by
$$
\phi_{\lambda}(z)=(1-\mid z\mid)^{\frac{i\lambda+n}{2}}\quad_{2}F_{1}(\frac{i\lambda+n}{2},\frac{i\lambda+n}{2},n;\mid z\mid^{2}).
$$
Using Pfaff's transformation \cite{Gr}
\begin{equation}
_{2}\digamma _{1}\left( a,b,c;x\right) =\left( 1-x\right) ^{-b}\, \, _{2}\digamma
_{1}\left( b,c-a,c;\frac{x}{x-1}\right)
\end{equation}
we rewrite $\phi_{-\lambda}$ as
\begin{equation}
 \phi_{-\lambda}(z)=\quad_{2}F_{1}\left(\frac{-i\lambda+n}{2},\frac{i\lambda+n}{2},n;\frac{\mid z\mid^{2}}{\mid z\mid^{2}-1}\right)
\end{equation}
So that returning back to (5.5) we get
\begin{equation}
\mathcal{F}[h_{m}^{\nu ,n}](\lambda)=2n\int\limits_{0}^{1}\frac{\rho ^{2n-1}}{\left( 1-\rho ^{2}\right)
^{n+1-2\left( \nu -m\right) }}\left( P_{m}^{\left( n-1,2(\nu -m)-n\right)
}\left( 1-2\rho ^{2}\right) \right) ^{2}
\end{equation}
\begin{equation*}
\times _{2}\digamma _{1}\left( \frac{n+i\lambda }{2},\frac{n-i\lambda }{2},n;%
\frac{\rho ^{2}}{\rho ^{2}-1}\right) d\rho
\end{equation*}
To calculate this last integral, we first use a linearisation of the
square of Jacobi polynomial in (5.8) by making appeal to the
following Clebsh-Gordon type formula see \cite[p. 611]{C},
\begin{equation}
P_{s}^{\left( \kappa ,\epsilon \right) }\left( u\right) P_{l}^{\left( \tau
,\eta \right) }\left( u\right) =\sum\limits_{k=0}^{s+l}A_{s,l}\left(
k\right) P_{k}^{\left( \alpha ,\delta \right) }\left( u\right)
\end{equation}
for the particular case of parameters $s=l=m$, $\kappa =\tau =\alpha =n-1$, $%
\epsilon =\eta =2\left( \nu -m\right) -n$ \, and $\delta = 2\left( \nu
-m\right) .$ In our setting, the linearisation coefficients $A_{s,l}\left(
k\right) $ are of the form
\begin{equation}
A_{k}^{\nu ,n,m}=\frac{\left( 2\left( \nu -m\right) +n\right) _{k}\left(
n\right) _{2m}\left( 2k+2\left( \nu -m\right) +n\right) \left( -1\right)
^{k}\left( 2m\right) !\left( \left( 2\left( \nu -m\right) \right)
_{2m}\right) ^{2}}{\left( n\right) _{k}\left( 2\left( \nu -m\right)
+n\right) _{2m+k+1}\left( m!\right) ^{2}\left( 2m-k\right) !\left( \left(
2\left( \nu -m\right) \right) _{m}\right) ^{2}}
\end{equation}
\begin{equation*}
\times \digamma _{2:1}^{2:2}\left(
\begin{array}{c}
-2m+k,-2\nu -k-n:-m,-n-m+1;-m,-m-n+1 \\
-2m,-2m-n+1:1-2\nu ,1-2\nu
\end{array}
\mid 1,1\right)
\end{equation*}
Here $\digamma _{l:l^{\prime }}^{p:p^{\prime }}\left( .\right) $ denotes the
Kamp\'{e} de F\'{e}riet double hypergeometric function defined by  \cite[p. 63]{S}
\begin{equation}
\digamma _{l:l^{\prime }}^{p:p^{\prime }}\left(
\begin{array}{c}
\left( a_{p}\right) :\left( b_{p^{\prime }}\right) ,\left( c_{p^{\prime
}}\right)  \\
\left( d_{l}\right) :\left( \kappa _{l^{\prime }}\right) ,\left( \varrho
_{l^{\prime }}\right)
\end{array}
\mid x,y\right) =\sum\limits_{q,s=0}^{+\infty }\frac{\left[ a_{p}\right]
_{q+s}\left[ b_{p^{\prime }}\right] _{q}\left[ c_{p^{\prime }}\right] _{s}}{%
\left[ d_{l}\right] _{q+s}\left[ \kappa _{l^{\prime }}\right] _{q}\left[
\varrho _{l^{\prime }}\right] _{s}}\frac{x^{q}}{q!}\frac{y^{s}}{s!}
\end{equation}
where $\left[ a_{p}\right] _{s}=\prod_{j=1}^{p}\left( a_{j}\right) _{s}$ in
which $\left( x\right) _{s}=x\left( x+1\right) ...\left( x+s-1\right) $ is
the Pochhammer symbol. Therefore, inserting
\begin{equation}
\left( P_{m}^{\left( n-1,2(\nu -m)-n\right) }\left( 1-2\rho ^{2}\right)
\right) ^{2}=\sum\limits_{k=0}^{2m}A_{k}^{\nu ,n,m}P_{k}^{\left( n-1,2(\nu
-m)\right) }\left( 1-2\rho ^{2}\right)
\end{equation}
into (5.8) the Fourier-Helgason transform  of $h_{m}^{\nu,n}$ takes the form
\begin{equation}
\mathcal{F}[h_{m}^{\nu ,n}](\lambda)
=\sum\limits_{k=0}^{2m}A_{k}^{\nu ,n,m}\frak{I}_{k}^{\nu ,n,m}\left( \lambda
\right),
\end{equation}
where the last term in this sum is defined by the integral
\begin{eqnarray}
\frak{I}_{k}^{\nu ,m}\left( \lambda \right)  &=&\int\limits_{0}^{1}\frac{%
2n\rho ^{2n-1}}{\left( 1-\rho ^{2}\right) ^{n+1-2\left( \nu -m\right) }}%
P_{k}^{\left( n-1,2(\nu -m)\right) }\left( 1-2\rho ^{2}\right)
 \\
&&\times _{2}F_{1}\left( \frac{1}{2}\left( n+i\lambda \right) ,\frac{1}{2}%
\left( n-i\lambda \right) ,n;\frac{\rho ^{2}}{\rho ^{2}-1}\right) d\rho
\notag
\end{eqnarray}
To calculate this last integral we make the change of variable $\rho =\tanh
t.$ Therefore (5.14) \ reads
\begin{eqnarray}
 \frak{I}_{k}^{\nu ,m}\left( \lambda \right)  &=&\int\limits_{0}^{+\infty
}2n\left( \sinh t\right) ^{2n-1}P_{k}^{\left( n-1,2(\nu -m)\right) }\left(
1-2\tanh ^{2}t\right)\\
&&\times \left( \cosh t\right) ^{-4\left( \nu -m\right) +1}._{2}F_{1}\left(
\frac{n+i\lambda }{2},\frac{n-i\lambda }{2},n;-\sinh ^{2}t\right) dt \notag
\end{eqnarray}
Now, we make use of the result established by Koornwinder \cite{K}
\begin{align}
& \int\limits_{0}^{+\infty }(\cosh t)^{-\alpha +\beta -\delta -\mu ^{\prime
}-1}\left( \sinh t\right) ^{2\alpha +1}P_{k}^{\left( \alpha ,\delta \right)
}\left( 1-2\tanh ^{2}t\right)  \nonumber \\
& \qquad \times _{2}F_{1}\left( \frac{\alpha +\beta +1+i\lambda }{2},\frac{\alpha
+\beta +1-i\lambda }{2},\alpha +1;-\sinh ^{2}t\right) dt  \notag
\\
&=\frac{\Gamma \left( \alpha +1\right) \left( -1\right) ^{k}\Gamma \left(
\frac{1}{2}\left( \delta +\mu ^{\prime }+1+i\lambda \right) \right) \Gamma
\left( \frac{1}{2}\left( \delta +\mu ^{\prime }+1-i\lambda \right) \right) }{%
k!\Gamma \left( \frac{1}{2}\left( \alpha +\beta +\delta +\mu ^{\prime
}+2\right) +k\right) \Gamma \left( \frac{1}{2}\left( \alpha -\beta +\delta
+\mu ^{\prime }+2\right) +k\right) }
\\
& \qquad
\times W_{k}\left( \frac{1}{4}\lambda ^{2};\frac{1}{2}\left( \delta +\mu
^{\prime }+1\right) ,\frac{1}{2}\left( \delta -\mu ^{\prime }+1\right) ,%
\frac{1}{2}\left( \alpha +\beta +1\right) ,\frac{1}{2}\left( \alpha -\beta
+1\right) \right) \nonumber
\end{align}
where $\beta ,\delta ,\lambda \in \mathbb{R}$, $\alpha ,\delta >-1$ , $%
\delta +\Re (\mu )^{\prime }>-1$ and  $W_{k}\left( .\right) $ is the
Wilson polynomial given in terms of the $_{4}F_{3}$-sum as (\cite{An},p. 158):
\begin{align}
W_{k}\left( x^{2},a,b,c,d\right) &:=\left( a+b\right) _{k}\left( a+c\right)
_{k}\left( a+d\right) _{k}  \notag
\\ & \qquad
\times _{4}F_{3}\left(
\begin{array}{c}
-k,k+a+b+c+d-1,a+ix,a-ix \\
a+b,a+c,a+d
\end{array}
\mid 1\right)
\end{align}
for the parameters $\alpha =n-1,\delta =2(\nu -m)-n,\beta =0$ and $\mu
^{\prime }=2(\nu -m)-n-1.$We find that
\begin{equation}
\frak{I}_{k}^{\nu ,m}\left( \lambda \right) =\frac{2n\Gamma \left( n\right)
\left( -1\right) ^{k}}{k!\Gamma ^{2}\left( 2\left( \nu -m\right) +k\right) }%
\left| \Gamma \left( 2\left( \nu -m\right) -\frac{n}{2}+i\frac{\lambda }{2}%
\right) \right| ^{2}
\end{equation}
\begin{equation*}
\times W_{k}\left( \frac{1}{4}\lambda ^{2};2\left( \nu -m\right) -\frac{n}{2}%
,1+\frac{n}{2},\frac{n}{2},\frac{n}{2}\right)
\end{equation*}
Summarizing the above calculations
\begin{equation}
\mathcal{F}\left[ h_{m}^{\nu ,n}\right] \left( \lambda \right) =\left|
\Gamma \left( 2\left( \nu -m\right) -\frac{n}{2}+i\frac{\lambda }{2}\right)
\right| ^{2}
\end{equation}
\begin{equation*}
\times \sum\limits_{k=0}^{2m}\gamma _{k}^{\nu ,n,m}W_{k}\left( \frac{\lambda
}{4}^{2};2\left( \nu -m\right) -\frac{n}{2},1+\frac{n}{2},\frac{n}{2},\frac{n%
}{2}\right)
\end{equation*}
with the constants
\begin{equation}
\gamma _{k}^{\nu ,n,m}:=\frac{2m!\Gamma \left( n\right) \left( 2\left( \nu
-m\right) -n\right) \Gamma \left( 2\nu -m\right) \left( -1\right)
^{k}A_{k}^{\nu ,n,m}}{\Gamma \left( n+m\right) \Gamma \left( 2\nu
-m-n+1\right) k!\Gamma ^{2}\left( 2\left( \nu -m\right) +k\right) },
\end{equation}
where the constants $A_{k}^{\nu ,n,m}$ is given by (5.10).
Finally, replacing $\lambda $ by $\sqrt{-\Delta _{\mathbb{B}^{n}}-n^{2}}$,
we arrive at the announced result.\\

\noindent \textbf{Remark 5.1.} Setting $m=0$ in the  formula (5.3) in Theorem 5.1, we recover the result of Peetre \cite{P}.\\

\noindent \textbf{Remark 5.2.} We should note that\textbf{\ }the transform\ $\frak{B}%
_{m}^{\nu ,n}$ \ have been expressed in \cite{Gh} as a function of
the Laplace-Beltrami operator $\Delta _{\mathbb{B}^{n}}$ in terms of the $%
_{3}\digamma _{2}$-sum as
\begin{align}
\frak{B}_{m}^{\nu ,n} &=\sum\limits_{j=0}^{2m}C_{j}^{\nu ,n,m}\frac{\Gamma
\left( 2\left( \nu -m\right) -\frac{1}{2}\left( n-i\sqrt{-\Delta _{\mathbb{B}%
^{n}}-n^{2}}\right) \right) }{\Gamma \left( 2\left( \nu -m\right) +j+\frac{1%
}{2}\left( n+i\sqrt{-\Delta _{\mathbb{B}^{n}}-n^{2}}\right) \right) }
\\
& \quad \times _{3}\digamma _{2}\left[
\begin{array}{c}
\frac{1}{2}\left( n+i\sqrt{-\Delta _{\mathbb{B}^{n}}-n^{2}}\right) ,n+j,%
\frac{1}{2}\left( n+i\sqrt{-\Delta _{\mathbb{B}^{n}}-n^{2}}\right) \\
\left( \nu -m\right) +j+\frac{1}{2}\left( n+i\sqrt{-\Delta _{\mathbb{B}%
^{n}}-n^{2}}\right) ,n
\end{array}
\mid 1\right] \nonumber
\end{align}
where
\begin{align}
C_{j}^{\nu ,n,m}&=\frac{\left( 2\left( \nu -m\right) -n\right) \Gamma \left(
n+m\right) \left( -1\right) ^{j}\Gamma \left( n+j\right) }{m!\Gamma \left(
2\nu -n-m+1\right) \Gamma \left( 2\nu -n\right) }
\\
&\times \sum\limits_{p=\max \left( 0,j-m\right) }^{\min \left( m,j\right) }%
\frac{\left( m!\right) ^{2}\Gamma \left( 2\nu -m\right) \Gamma \left( 2\nu
-m+j-p\right) }{\left( j-p\right) !\left( m+p-j\right) !p!\left( m-p\right)
!\Gamma \left( n+j-p\right) \Gamma \left( n+p\right) }. \nonumber
\end{align}

\end{document}